\documentclass[12pt]{article}
\usepackage[cp866]{inputenc}
\usepackage[psamsfonts]{amssymb}
\usepackage{amsmath,amsfonts,amssymb,amscd}
\usepackage{graphicx}
\newtheorem{theorem}{Theorem}[section]

\newtheorem{lemma}[theorem]{Lemma}
\newtheorem{prop}[theorem]{Proposition}
\newtheorem{problem}[theorem]{Problem}
\newtheorem{remark}[theorem]{Remark}

\begin{document}

\title{Groups of virtual trefoil and Kishino knots}

\author{{V.~G.~Bardakov, Yu. A. Mikhalchishina, M.~V.~Neshchadim }}%

\maketitle {\small}

\vspace{0.8cm}

\begin{abstract}
In the paper \cite{Mih} for an arbitrary virtual link $L$  three groups $G_{1,r}(L)$, $r>0$, $G_{2}(L)$ and $G_{3}(L)$ were defined.
In the present paper  these groups for the virtual trefoil are investigated. The structure of these groups are found out and the fact that some of them are not isomorphic to each other is proved. Also we prove that $G_3$ distinguishes the Kishino knot from the trivial knot. The fact that these groups have the lower central series which does not stabilize on the second term is noted. Hence we have a possibility to study these groups using quotients by terms of the lower central series and to construct representations of these groups in rings of formal power series. It  allows to construct an invariants for virtual knots.
\end{abstract}


\section{Introduction}

\vspace{0.5cm}

Unlike classical knots virtual ones are not topological objects therefore there is no natural definition of a virtual knot group. There exist several approaches
to the definition (see \cite{Ka, SW, Bar-1, CSW, BB, BDG, BMN, BMN-1, Mih}). Some of the approaches use a knot diagram, other ones use representations of the
virtual braid group by automorphisms of some groups. The following natural problems occur of how to compare these groups, how to describe isomorphic groups and
how to construct homomorphisms of some groups into another ones. A basis of such an analysis was put in \cite{BMN-1}.

In the present work groups defined by Yu. A. Mikhalchishina \cite{Mih} are analyzed. In \cite{Mih} representations $W_{1,r}$, $r > 0$, $W_2$ and $W_3$ of the virtual
braid group $VB_n$ into the automorphism group $\mathrm{Aut}(F_{n+1})$ of a free group $F_{n+1} = \langle y, x_1, x_2, \ldots, x_n \rangle$ of rank $n+1$ were
constructed. These representations extend Wada representations \cite{W} $w_{1,r}$, $r > 0$, $w_2$, $w_3$  of the braid group $B_n$ into the automorphism group
$\mathrm{Aut}(F_{n})$ of a free group $F_{n} = \langle  x_1, x_2, \ldots, x_n \rangle$ of rank $n$. Using these representations, for each virtual link $L$ the
groups $G_{1,r}(L)$, $G_{2}(L)$ and $G_{3}(L)$ were defined. In particular, generators and defining relations of the groups for the virtual trefoil were
written. In the present work the fact that groups of the virtual trefoil $T_v$ $G_{1,r}(T_v)$ and $G_{2}(T_v)$ are not isomorphic to each other will be proven.

Note that the most of virtual knot groups do not distinguish the group of the Kishino knot from the trivial one. The group defined in \cite{CSW} is an
exception. In the paper \cite{Mih} the fact that groups constructed using representations $W_{1,r}$ and  $W_{2}$ do not distinguish the Kishino knot from the
trivial one is stated. And the natural question whether the group constructed using the representation $W_{3}$ distinguishes the Kishino knot from the trivial
one  was formulated in \cite{Mih}. In present work the positive answer on the question is given.

As it is well known if  $G$ is a group of a classical knot, then its commutator subgroup coincides with the third term of lower central series, i.~e. $[G,G] =
[[G,G], G]$. Therefore, the factorization by terms of the lower central series cannot be used for distinguish of knots. On the other hand as it was pointed out in
the paper \cite{BM} some of virtual link groups are residually nilpotent groups. Recall that a group $G$ is referred to as \textit{residually  nilpotent} if for
every nontrivial element $g$ from $G$ there exists a homomorphism $\varphi : G \longrightarrow N$ on a nilpotent group $N$ such that $\varphi(g) \not= 1$. It is
easy to observe that $G$ is residually nilpotent if and only if the intersection of all terms of its lower central series is a trivial group.

In the paper \cite{BN} was noted that the virtual trefoil group from \cite{Bar-1} has different first five terms of lower central series. This fact allows to
construct new invariants of virtual links. The question occurs whether virtual link groups are residually nilpotent or not. Moreover in \cite{BN}
representations of some virtual knot groups into finite algebras were constructed.  In the
present paper this approach is being developed.

Let us describe the paper content. In Section 2 we recall representations constructed in \cite{Mih} of the virtual braid group $VB_n$ into the group of automorphisms of free groups which are extensions of the Wada representations of braid group $B_n$. Also we give presentations of three groups $G_{1,r}$, $G_2$ and $G_3$ of the virtual trefoil. In Section 3  the structure of these groups is found out. In Section 4, using quotients by the second commutator subgroups, we prove that $G_{1,r}$ is not isomorphic to $G_2$; we study the quotients by the terms of the lower central series and prove that $\gamma_2 (G_2) / \gamma_3 (G_2)$ is the cyclic group of order 2 and the quotient $\gamma_k (G_2) / \gamma_{k+1} (G_2)$ has exponent 4 for all $k \geq 2$. Studying the lower central series for $G_{1,r}$, we prove that $G_{1,r} / \gamma_{4} (G_{1,r}) \cong F_2 / \gamma_{4} (F_2)$ and $\gamma_4 (G_{1,r}) / \gamma_{5} (G_{1,r}) \cong \mathbb{Z} \times \mathbb{Z} \times \mathbb{Z}_r$. In Section 5 we prove that $G_3(Ki)$ of the Kishino knot $Ki$ is not isomorphic to $G_3(U)$ of the trivial knot $U$. In Section 6 we construct a representation of an arbitrary finitely presented group into a ring of formal power series with non-commutative variables and prove that the image does not depend on the presentation of group; we construct representations of $G_{1,r}$ and $G_2$ into some finite dimensional algebra.

\textbf{Acknowledgement}. This work was supported by  the Russian Science Foundation (project No. 16-41-02006). The authors are also grateful to E. I. Timoshenko for the consultation on metabelian groups.

\vspace{0.8cm}

\section{Extensions of Wada representations and virtual trefoil groups}

\vspace{0.5cm}

In the paper \cite{Mih}  the following three representations of the virtual braid group $VB_n$ into the automorphism group
$\mathrm{Aut}(F_{n+1})$ of the rank $n+1$ free group $F_{n+1} = \langle y, x_1, x_2, \ldots, x_n \rangle$ were constructed.

 1. The representation $W_{1,r},\, r > 0$ is defined by the action on the generators
$$
W_{1,r}(\sigma_{i}) :\left\{ \begin{array}{ll}
    x_{i} \longmapsto x_{i}^{r} x_{i+1} x_{i}^{-r}, \\
    x_{i+1} \longmapsto x_{i},\\
    \end{array}
    \right.
~~~ W_{1,r}(\rho_i) :\left\{ \begin{array}{ll}
    x_{i} \longmapsto  x_{i+1}^{y^{-1}}, \\
    x_{i+1} \longmapsto x_{i}^{y}. \\
        \end{array}
 \right.$$ Here and onward we point out only nontrivial actions on generators assuming that other generators are fixed.

2. The representation $W_{2}$ is defined by the action on the generators $$W_{2}(\sigma_{i}) :\left\{ \begin{array}{ll}
    x_{i} \longmapsto x_{i} x_{i+1}^{-1} x_{i}, \\
    x_{i+1} \longmapsto x_{i},\\
  \end{array}
\right. ~~~ W_{2}(\rho_i) :\left\{ \begin{array}{ll}
    x_{i} \longmapsto  x_{i+1}^{y^{-1}}, \\
    x_{i+1} \longmapsto x_{i}^{y}. \\
      \end{array}
 \right. $$

3. The representation $W_{3}$ is defined by the action on the generators $$W_{3}(\sigma_{i}) :\left\{ \begin{array}{ll}
    x_{i} \longmapsto x_{i}^{2} x_{i+1}, \\
    x_{i+1} \longmapsto x_{i+1}^{-1} x_{i}^{-1} x_{i+1}, \\
    \end{array}
\right. ~~~ W_{3}(\rho_i) :\left\{ \begin{array}{ll}
    x_{i} \longmapsto  x_{i+1}^{y^{-1}}, \\
    x_{i+1} \longmapsto x_{i}^{y}. \\
        \end{array}
 \right. $$
These representations extend Wada representations \cite{W} $w_{1,r}$, $r > 0$, $w_2$, $w_3$ of the braid
group $B_n$ into the automorphisms group $\mathrm{Aut}(F_{n})$ of the rank $n$ free group $F_{n} = \langle  x_1, x_2, \ldots, x_n \rangle$.

In the paper \cite{Mih} for each virtual link three types of groups were defined: $G_{1,r}(L)$, $G_{2}(L)$ and $G_{3}(L)$ that correspond to described
representations. And the fact that these groups are invariants of a virtual link $L$ is proved. In particular,
for the virtual trefoil these groups have the form:

 \smallskip

$ G_{1,r} = \langle\, x,\,y\, \| \, x^{-r}y^{-1}xyx^{r}=y^{-2}x^{r}yxy^{-1}x^{-r}y^2 \,  \rangle,~~r > 0, $ \smallskip

$ G_2 = \langle\, x,\,y\, \| \, xy^{-1}x^{-1}yx=y^{-2}xyx^{-1}y^{-1}xy^2 \,  \rangle, $ \smallskip

$ G_3 = \langle\, x_1,\,x_2,\,y\, \| \, x_1=yx_2^{-1}x_1^{-1}x_2^{-1}x_1^{-1}x_2y^{-1},\,
                                    x_2=y^{-1}x_1^{2}x_2x_1x_2y \,  \rangle.
$\\
Note that the group $G_{1,1}$ is isomorphic to the group found out in \cite{BB} (see also \cite{BMN}).


\vspace{0.8cm}

\section{Structure of the virtual trefoil groups}

\vspace{0.5cm}

The structure  of the virtual trefoil groups gives

\begin{theorem} 1) The group $G_{1,r}$ is the semi-direct product $H \leftthreetimes \mathbb{Z}$, where
$$
H =\ldots \underset{B_{-2}}{\ast} A_{-1}\underset{B_{-1}}{\ast} A_{0}
  \underset{B_{0}}{\ast} A_{1} \underset{B_{1}}{\ast} \ldots.
$$
is the amalgamated free product with
$$
 A_i=\langle\, x_{i-1},\,x_i,\, x_{i+1} \| \,
      [x_i, x_{i-1}^{\varepsilon} x_{i+1}^{-\varepsilon}]=1 \,  \rangle
 \cong (\langle\, x_i\,\rangle \times
       \langle \,x_{i-1}^{\varepsilon} x_{i+1}^{-\varepsilon} \,\rangle )\ast\langle\,x_{i+1}\,\rangle
    = (\mathbb{Z} \times \mathbb{Z})\ast\mathbb{Z} ~\mbox{if}~r=\varepsilon=\pm 1,
$$
and
$$ A_i=\langle\, x_{i-1},\,x_i,\, x_{i+1} \| \,
      x_{i}^{-1} (x_{i-1}^{r}x_{i+1}^{-r})  x_{i}=x_{i-1}^{r}x_{i+1}^{-r}\,  \rangle ~\mbox{if}~r=\varepsilon=\pm 1;
$$
all $B_i$  are isomorphic to $F_2$.

2) $G_{2}$ is the semi-direct product $H \leftthreetimes \mathbb{Z}$, where
$$
H =\ldots \underset{B_{-2}}{\ast} A_{-1}\underset{B_{-1}}{\ast} A_{0}
  \underset{B_{0}}{\ast} A_{1} \underset{B_{1}}{\ast} \ldots.
$$
is the amalgamated free product with
$$
A_i=\langle\, x_{i-1},\,x_i,\, x_{i+1} \| \,
     (x_{i-1} x_{i+1}^{-1})^{x_{i} x_{i-1}^{-1}} = (x_{i-1} x_{i+1}^{-1})^{-1}\,  \rangle
$$
and  all $B_i$  are isomorphic to $F_2$.

3) $G_3$ is an HNN-extension:
$$
G_3 = \langle x_1, x_2, y ~||~ y^{-1} A y = B, \varphi \rangle
$$
with the base group $\langle x_1, x_2 \rangle = F_2$ and associated subgroups $A = \langle x_2 x_1 x_2, x_1 \rangle$, $B = \langle x_1 x_2 x_1, x_2 \rangle$, which are also isomorphic to $F_2$.

\end{theorem}

\textbf{Proof.}  1) Consider the group $G_{1,r}$. Let $$ x_i=y^{-i}xy^{i},\quad i \in \mathbb{Z}. $$ Then the defining relation of the group $G_{1,r}$ can be
rewritten in the form $$ x_0^{-r}x_1x_0^{r}=x_2^{-r}x_1x_2^{r}. $$ Conjugating the relation by  $y^{i}$, $i \in \mathbb{Z}$, we have $$
x_i^{-r}x_{i+1}x_i^{r}=x_{i+2}^{-r}x_{i+1}x_{i+2}^{r},\quad i \in \mathbb{Z}. $$ Let us define the following groups $$ H =\langle\, x_i,\,\, i \in \mathbb{Z} \|
\,
      x_i^{-r}x_{i+1}x_i^{r}=x_{i+2}^{-r}x_{i+1}x_{i+2}^{r},\,\, i \in \mathbb{Z} \,  \rangle,
$$ $$ A_i=\langle\, x_{i-1},\,x_i,\, x_{i+1} \| \,
      x_{i-1}^{-r}x_{i}x_{i-1}^{r}=x_{i+1}^{-r}x_{i}x_{i+1}^{r} \,  \rangle,
$$ $$ B_{i}=\langle\, x_i,\, x_{i+1}  \,  \rangle. $$ Then the group $G_{1,r}$ can be decomposed into the semidirect product $$ G_{1,r} = H \leftthreetimes
\mathbb{Z}, $$ where $\mathbb{Z}=\langle y \rangle$ is the infinite cyclic group generated by the element $y$, the group $H$ is the amalgamated free product $$
H =\ldots \underset{B_{-2}}{\ast} A_{-1}\underset{B_{-1}}{\ast} A_{0}
  \underset{B_{0}}{\ast} A_{1} \underset{B_{1}}{\ast} \ldots.
$$ Since $A_i$ is a group with one relation the subgroup $B_{i}$ is a free group of rank 2.

The defining relation of the subgroup $A_i$ of the group $G_{1,r}$ can be rewritten in the form
$$ [x_i, x_{i-1}^{r}x_{i+1}^{-r}]=1. $$ If
$r=\varepsilon=\pm 1$, then
$$
 A_i=\langle\, x_{i-1},\,x_i,\, x_{i+1} \| \,
      [x_i, x_{i-1}^{\varepsilon} x_{i+1}^{-\varepsilon}]=1 \,  \rangle
 \cong (\langle\, x_i\,\rangle \times
       \langle \,x_{i-1}^{\varepsilon} x_{i+1}^{-\varepsilon} \,\rangle )\ast\langle\,x_{i+1}\,\rangle
    = (\mathbb{Z} \times \mathbb{Z})\ast\mathbb{Z};
$$
if $r\neq \pm 1$, then the group $A_i$ is the $HNN$-extension: $$ A_i=\langle\, x_{i-1},\,x_i,\, x_{i+1} \| \,
      x_{i}^{-1} (x_{i-1}^{r}x_{i+1}^{-r})  x_{i}=x_{i-1}^{r}x_{i+1}^{-r}\,  \rangle.
$$

2) Consider the group $G_{2}$. Let $$ x_i=y^{-i}xy^{i},\quad i \in \mathbb{Z}. $$ Then the defining relation of the group $G_2$ can be rewritten in the form $$
x_0x_1^{-1}x_0=x_2x_1^{-1}x_2. $$ Conjugating this relation by $y^{i}$, $i \in \mathbb{Z}$, we obtain $$ x_i x_{i+1}^{-1}x_i=x_{i+2}x_{i+1}^{-1}x_{i+2},\quad i
\in \mathbb{Z}. $$ Define the following groups $$ H = \langle\, x_i,\,\, i \in \mathbb{Z} \| \,
      x_i x_{i+1}^{-1}x_i=x_{i+2}x_{i+1}^{-1}x_{i+2},\,\, i \in \mathbb{Z} \,  \rangle,
$$ $$ A_i=\langle\, x_{i-1},\,x_i,\, x_{i+1} \| \,
      x_{i-1}x_{i}^{-1}x_{i-1}=x_{i+1}x_{i}^{-1}x_{i+1} \,  \rangle,
$$ $$ B_{i}=\langle\, x_i,\, x_{i+1}  \,  \rangle. $$ Then the group $G_2$ can be decomposed into the semidirect product $$ G_2=H \leftthreetimes \mathbb{Z}, $$
where $\mathbb{Z}=\langle y \rangle$ is the infinite cyclic group generated by the element $y$, the group $H$ is the amalgamated free product $$ H = \ldots
\underset{B_{-2}}{\ast} A_{-1}\underset{B_{-1}}{\ast} A_{0}
  \underset{B_{0}}{\ast} A_{1} \underset{B_{1}}{\ast} \ldots.
$$
The defining relation of the subgroup $A_i$ of the group $G_2$  can be rewritten in the form $$ x_i^{-1} (x_{i-1}x_{i+1}^{-1})x_i=x_{i-1}^{-1}x_{i+1}
\Leftrightarrow (x_{i-1} x_{i+1}^{-1})^{x_{i} x_{i-1}^{-1}} = (x_{i-1} x_{i+1}^{-1})^{-1}. $$ Therefore, the group $A_i$ is the $HNN$-extension and it contains
the subgroup which is isomorphic to the fundamental group of the Klein bottle.

Since $A_i$ is a group with one relation the subgroup $B_{i}$ is a free group of rank 2.

3) Rewrite relations of the group $G_3$ in the form $$ x_1^y=x_2^{-1}x_1^{-1}x_2^{-1}x_1^{-1}x_2,\quad
   (x_1^{2}x_2x_1x_2)^{y}=x_2.
$$ Let $$ A=\langle\, x_1^{2}x_2x_1x_2,\, x_1 \,  \rangle=\langle\, x_2x_1x_2,\, x_1 \,  \rangle, $$ $$
 B=\langle\, x_2^{-1}x_1^{-1}x_2^{-1}x_1^{-1}x_2,\, x_2 \,  \rangle=
 \langle\, x_1x_2x_1,\, x_2 \,  \rangle.
$$ Then the group $G_3$ is the $HNN$-extension with the stable letter $y$. Groups $A$ and $B$ are 2-generated free subgroups in the free group
$\langle\, x_1,\, x_2 \,  \rangle$ of rank 2 and the isomorphism $\varphi : A \rightarrow B$ is defined by the rule
$$
x_1 \mapsto x_2^{-1} x_1^{-1} x_2^{-1} x_1^{-1} x_2,~~~ x_1^{2} x_2 x_1  x_2 \mapsto x_2.
$$

\hfill $\square$ \medskip

\medskip

A natural question whether groups $G_{1,r}, G_2$ and $G_3$ are isomorphic to each other occurs. The next section is dedicated to the answer on this question.

\vspace{0.8cm}

\section{Isomorphism problem}

\vspace{0.5cm}

In the present section the fact that some of virtual trefoil groups listed above are not isomorphic to each other will be proven. For that purpose different
approaches will be used: a quotient by the second commutator subgroup, a quotient by terms of lower central series.

\subsection{Quotient group by the second commutator subgroup}  1) Find the quotient group $G_2/G_2''$. In order to do that note the relation of the group $G_2$: $$
xy^{-1}x^{-1}yx=y^{-2}xyx^{-1}y^{-1}xy^2 $$ is equivalent to the following relation $$ [x,y]^{2}[x,y]^{2y}=1. $$ Therefore, the quotient group $G_2/G_2''$ can
be considered as a cyclic module $G_2'/G_2''=\langle [x,y]G_2'' \rangle$ over the Laurent  ring $\mathbb{Z}[x^{\pm 1}, y^{\pm 1}]$. In this module relation can be written in
the form $$ [x,y]\cdot 2(1+y)=0. $$

2) Find the quotient group $G_{1,r} / G_{1,r}''$. The following commutator identities will be used $$ [ab,c]=[a,c]^b[b,c],~~~[a,bc]=[a,c][a,b]^c,~~~
[a,b^{-1}]=[b,a]^{b^{-1}},~~~[a,b]^{-1}=[b,a]. $$

Using these identities, the relation of the group $G_{1,r}$ is transformed as follows: $$ x^{-r}y^{-1}xyx^{r}=y^{-2}x^{r}yxy^{-1}x^{-r}y^2 \Leftrightarrow $$ $$
\Leftrightarrow x^{-1}(yx^{r})^{-1}x(yx^{r})=x^{-1}(y^{-1}x^{-r}y^2)x (y^{-1}x^{-r}y^2) \Leftrightarrow $$ $$\Leftrightarrow [x,yx^{r}]=[x, y^{-1}x^{-r}y^2] \Leftrightarrow $$ $$\Leftrightarrow [x,y]^{x^{r}}=[x, y^2][x,
y^{-1}x^{-r}]^{y^2}. $$ Transform the right part of the last relation: $$ [x, y^2][x, y^{-1}x^{-r}]^{y^2}=[x, y][x, y]^y[x, y^{-1}]^{x^{-r}y^2}= $$ $$ =[x, y][x,
y]^y[y,x ]^{y^{-1}x^{-r}y^2}=[x, y][x, y]^y[x,y ]^{-y^{-1}x^{-r}y^2}. $$ Thus the relation of the group $G_{1,r}$ by modulo of the second commutator subgroup
$G_{1,r}''$ has the form $$ [x,y]\cdot (1+y-yx^{-r}-x^{r})=0 $$ or $$ [x,y]\cdot (1-x^{-r})(y-x^{r})=0. $$

\vspace{0.5cm}

To prove that the group $G_{2}$ is not isomorphic to any group $G_{1,r}$ the following lemma is required.

\begin{lemma} If the automorphism of the free metabelian group $F_2/F_2''$ is defined by the action on the  generators: $$ x \mapsto
x^{\alpha}y^{\beta}[x,y]^{\gamma},\quad y \mapsto x^{a}y^{b}[x,y]^{c}, ~ \mbox{where}~\alpha, \beta, \gamma, a, b, c \in \mathbb{Z}, ~~\alpha b-\beta a=\pm 1, $$
then the commutator $z=[x,y] \in F_2/F_2''$ is transformed by the rule $$ z \mapsto z \cdot \left( c(1-x^{\alpha}y^{\beta})-\gamma (1-x^{a}y^{b}) + y^{\beta}
\frac{(1-x^{\alpha})}{(1-x)} \frac{(1-y^{b})}{(1-y)}- y^{b} \frac{(1-x^{a})}{(1-x)} \frac{(1-y^{\beta})}{(1-y)} \right). $$

\end{lemma}

\textbf{Proof.} All further computations are performed in the quotient group $F_2/F_2''$. We have $$ z \mapsto [x^{\alpha}y^{\beta}z^{\gamma},
x^{a}y^{b}z^{c}]= [x^{\alpha}y^{\beta}z^{\gamma},  z^{c}] [x^{\alpha}y^{\beta}z^{\gamma},  x^{a}y^{b}]= $$ $$ =[x^{\alpha}y^{\beta},z^{c}]
[x^{\alpha}y^{\beta}z^{\gamma},y^{b}][x^{\alpha}y^{\beta}z^{\gamma},x^{a}]^{y^{b}}= $$ $$ =[x^{\alpha},z^{c}]^{y^{\beta}} [y^{\beta},z^{c}]
[x^{\alpha}y^{\beta},y^{b}][z^{\gamma},y^{b}] [y^{\beta}z^{\gamma},x^{a}]^{y^{b}}= $$ $$ =[x^{\alpha},z^{c}]^{y^{\beta}}
[y^{\beta},z^{c}][x^{\alpha},y^{b}]^{y^{\beta}} [z^{\gamma},y^{b}][y^{\beta},x^{a}]^{y^{b}}[z^{\gamma},x^{a}]^{y^{b}}. $$ To transform this expression  the
following equations are used $$ [x^{p},y^{q}]=z\cdot \left( \frac{(1-x^{p})}{(1-x)} \frac{(1-y^{q})}{(1-y)}\right),\quad
 [x^{p},z^{q}]=z \cdot q(1-x^{p}),
$$ which are true for arbitrary integers  $p$, $q$. We obtain $$z \mapsto z\cdot c(1-x^{\alpha})y^{\beta}+ z\cdot c(1-y^{\beta})+ z\cdot y^{\beta}
\frac{(1-x^{\alpha})}{(1-x)} \frac{(1-y^{b})}{(1-y)}+ $$ $$ +z\cdot \gamma (y^{b}-1)+ z\cdot \gamma (x^{a}-1)y^{b}- z\cdot y^{b} \frac{(1-x^{a})}{(1-x)}
\frac{(1-y^{\beta})}{(1-y)}= $$ $$ = z \cdot \left( c(1-x^{\alpha}y^{\beta})-\gamma (1-x^{a}y^{b}) + y^{\beta} \frac{(1-x^{\alpha})}{(1-x)} \frac{(1-y^{b})}{(1-y)}- y^{b}
\frac{(1-x^{a})}{(1-x)} \frac{(1-y^{\beta})}{(1-y)} \right). $$ \hfill $\square$

Using this lemma, the following proposition is proved.

\begin{prop} The group $G_2$ is not isomorphic to any group $G_{1,r}$. \end{prop}

\textbf{Proof.} Using the transformation $$ x \mapsto x^{\alpha}y^{\beta}z^{\gamma},\quad y \mapsto x^{a}y^{b}z^{c} $$ to the expression $z\cdot 2(1+y)$, we
obtain $$ z\cdot 2(1+x^{a}y^{b}) \left( c(1-x^{\alpha}y^{\beta})-\gamma (1-x^{a}y^{b}) + y^{\beta} \frac{(1-x^{\alpha})}{(1-x)} \frac{(1-y^{b})}{(1-y)}- y^{b}
\frac{(1-x^{a})}{(1-x)} \frac{(1-y^{\beta})}{(1-y)} \right). $$ If $G_{1,r}$ is isomorphic to the group  $G_2$ the following equation should be true $$
2(1+x^{a}y^{b}) \left( c(1-x^{\alpha}y^{\beta})-\gamma (1-x^{a}y^{b}) + y^{\beta} \frac{(1-x^{\alpha})}{(1-x)} \frac{(1-y^{b})}{(1-y)}- y^{b} \frac{(1-x^{a})}{(1-x)}
\frac{(1-y^{\beta})}{(1-y)} \right)= $$ $$ =(1-x^{-r})(y-x^{r}) $$ for some parameters $\alpha$, $\beta$, $\gamma$, $a$, $b$, $c \in \mathbb{Z}$ and $\alpha b-\beta
a=\pm 1$. Setting $x=y=1$ and taking into account that $(1-x^{p}) / (1-x)=p$ for $x=1$ and for arbitrary integer $p$, we obtain the contradictory equation  $\pm
4=0$. \hfill $\square$

\medskip

\subsection{Lower central series} As it is well known for the classical knot group its commutator subgroup coincides with the third term of low central series. Let
us show that it is not the case for virtual knots. For example for the group $G_2$ it is true

\smallskip

\begin{prop} 1) The quotient group $\gamma_2 G_2/ \gamma_3 G_2$ is isomorphic to the cyclic group of order~4.

2) The quotient group $\gamma_k G_2/ \gamma_{k+1} G_2$ is a finite group of exponent 4 for an arbitrary integer $k \geq 2$.

3) The relation of the group $G_2$ by modulo $\gamma_4 G_2$ has the form $$ [y,x]^4=[x,y,y]^2. $$ \end{prop}

 \smallskip

\textbf{Proof.} 1) Transform the relation of the group $G_2$ $$ xy^{-1}x^{-1}yx=y^{-2}xyx^{-1}y^{-1}xy^2 \Rightarrow  x[y,x]=y^{-2}x[y^{-1},x]y^2
\Rightarrow  x[y,x]=xy^{-2}[y^{-2},x][y^{-1},x]y^2. $$ Therefore, the relation is transformed to the form $$ [y,x]=y^{-2}[y^{-2},x][y^{-1},x]y^2. $$ By
modulo $\gamma_3 G_2$ we have $$ [y,x]^4=1. $$

2) Since the quotient group $\gamma_k G_2 / \gamma_{k+1} G_2$ is generated by left-normalized commutators of the form $$ [\ldots [g_1, g_2], \ldots, g_k], ~~g_i \in
G_2, $$ raising to the forth power and using relations of the quotient group  $\gamma_k G_2 / \gamma_{k+1} G_2$, we obtain $$ [\ldots [g_1, g_2], \ldots, g_k]^4
= [\ldots [g_1, g_2]^4, \ldots, g_k] = 1 $$ due to the proved item 1).

3) Transform the right part of the relation $$ [y,x]=y^{-2}[y^{-2},x][y^{-1},x]y^2 $$ by  modulo $\gamma_4 G_2$. We have $$ y^{-2}[y^{-2},x][y^{-1},x]y^2=
[y^{-2},x^{y^2}][y^{-1},x^{y^2}]= [y^{-2},x[x,y^2]][y^{-1},x[x,y^2]]= $$ $$ =[y^{-2},[x,y^2]][y^{-2},x] [y^{-1},[x,y^2]][y^{-1},x]= [y,[x,y]]^{-4}[y^{-2},x]
[y,[x,y]]^{-2}[y^{-1},x]= $$ $$ =[y^{-2},x][y^{-1},x][x,y,y]^{6}= [y^{-1},x]^{y^{-1}}[y^{-1},x]^2[x,y,y]^{6}= $$ $$
=[y^{-1},x][y^{-1},x,y^{-1}][y^{-1},x]^2[x,y,y]^{6}= [y^{-1},x]^3[x,y,y]^{5}= $$ $$ =\left([x,y]^{y^{-1}}  \right)^3[x,y,y]^{5}= \left( [x,y] [x,y,y^{-1}]
\right)^3[x,y,y]^{5}= [x,y]^3[x,y,y]^{2}. $$ Thus the relation of the group $ G_2$ by modulo $\gamma_4 G_2$ has the form $$ [y,x]= [x,y]^3[x,y,y]^{2}, $$
which is equivalent to $$ [y,x]^4=[x,y,y]^2. $$ \hfill $\square$

\vspace{0.5cm}

\textbf{Remark.} The relation $$ [y,x]^4=[x,y,y]^2 $$ of the group  $ G_2$ by modulo $\gamma_4 G_2$ can be presented in the form $$
\left([y,x]^2\right)^y=[y,x]^{-2}. $$

Indeed $$ [y,x]^4=[x,y,y]^2=[[y,x]^{-2},y]=[y,x]^{2} y^{-1} [y,x]^{-2} y. $$

\bigskip

Consider terms of the lower central series of the group $G_{1,r}$. It is true

\begin{theorem} The first five terms of low central series of the group $ G_{1,r}$ are different from each other. Moreover $ G_{1,r}/  \gamma_4 G_{1,r} \cong
F_2/  \gamma_4 F_2$, $\gamma_4 G_{1,r}/  \gamma_5 G_{1,r} \cong \mathbb{Z} \times\mathbb{Z} \times\mathbb{Z}_r$. \end{theorem}

\textbf{Proof.} At first let us show that the relation of the group $G_{1,r}$ by modulo of the subgroup $\gamma_5 G_{1,r}$ has the form $$ [y,x,y,x]^r=
[y,x,x,y]^{r^2}. $$

For that purpose rewrite the relation of the group $G_{1,r}$ in the form $$ x^{yx^{r}}=x^{y^{-1}x^{-r}y^2}. $$ It is equivalent to  $$ [x,yx^{r}  y^{-2}x^{r}y]=1. $$ Using
commutator identities, transform the left side of this relation: $$ [x,yx^{r}  y^{-2}x^{r}y]=[x,(yx^{r}y^{-1}) (y^{-1}x^{r}y)]= [x,[y^{-1},x^{-r} ]x^{2r}
[x^{r},y]]= $$ $$ =[x, [x^{r},y]]  [x,[y^{-1},x^{-r} ]x^{2r} ]^{[x^{r},y]}= [x, [x^{r},y]]  [x,[y^{-1},x^{-r} ] ]^{x^{2r}[x^{r},y]}= $$ $$ =[x, [x^{r},y]]
[x,[y,x^{r} ]^{y^{-1}x^{-r}} ]^{x^{2r}[x^{r},y]}= [x, [x^{r},y]]  [x,[y,x^{r} ] [y,x^{r}, y^{-1}x^{-r}] ]^{x^{2r}[x^{r},y]}= $$ $$ =[x, [x^{r},y]]  [x,[y,x^{r},
y^{-1}x^{-r}] ]^{x^{2r}[x^{r},y]}
 [x,[y,x^{r} ]  ]^{[y,x^{r}, y^{-1}x^{-r}]x^{2r}[x^{r},y]}=
$$ $$ =[x, [x^{r},y]]  [x,[y,x^{r}, y^{-1}x^{-r}] ] [[x,[y,x^{r}, y^{-1}x^{-r}] ],x^{2r}[x^{r},y]]\times $$ $$
 \times
 [x,[y,x^{r} ]]  [[x,[y,x^{r} ]  ],[y,x^{r}, y^{-1}x^{-r}]x^{2r}[x^{r},y]]=
$$ (by modulo $\gamma_5 G_{1,r}$ we obtain) $$ =[x, [x^{r},y]]  [x,[y,x^{r}, y^{-1}x^{-r}] ]
 [x,[y,x^{r} ]]  [[x,[y,x^{r} ]  ],x^{2r}]=
$$ $$ =[x, [x^{r},y]][x,[y,x^{r} ]] [x,[y,x, y] ]^{-r}[x,[y,x, x] ]^{-r^2} [y,x,x,x]^{-2r^2}= $$ $$ =[x, [x^{r},y]] [x,[x^{r},y ]^{-1}] [y,x,y,x]^{r}
[y,x,x,x]^{-r^2}= $$ $$ =[x,[x^{r},y]]
   [[x^{r},y ],x]^{[x^{r},y ]^{-1}}
      [y,x,y,x]^{r}
          [y,x,x,x]^{-r^2}=
$$ $$ =[y,x,y,x]^{r} [y,x,x,x]^{-r^2}. $$ Recall that $[y,x,y,x]\equiv [y,x,x,y] (\mathrm{mod} \gamma_5 F_2)$. Indeed for $a=[y,x]$, $b=y$, $c=x$ from Hall
identity for  $ F_3/ \gamma_4 F_3$ we have $$ [a,b,c] [b,c,a] [c,a,b]=1. $$ Due to $a=[b,c]$ we obtain $$ [y,x,y,x] [x,[y,x],x]\equiv 1 (
\mathrm{mod} \gamma_5 G_{1,r}) $$ or $$ [y,x,y,x]\equiv [y,x,x,y]  ( \mathrm{mod} \gamma_5 F_2). $$ Therefore, the relation of the group $ G_{1,r}$ by modulo of
the subgroup $  \gamma_5 G_{1,r}$ has the form $$ [y,x,y,x]^r= [y,x,x,x]^{r^2}. $$

Basic commutators of weight 4 are of the form $c_1=[y,x,x,x]$,  $c_2=[y,x,x,y]$, $c_3=[y,x,y,y]$. Hence the relation of the quotient group $ \gamma_4 G_{1,r}/
\gamma_5 G_{1,r} $ can be written in the form $(c_2c_1^{-r})^r=1$. \hfill $\square$

\vspace{0.8cm}
\section{The proof of
non-triviality of Kishino knot}

\vspace{0.5cm}

It was proved in \cite{Mih} that groups $G_{1,r}(Ki)$ and  $G_{2}(Ki)$ cannot distinguish the Kishino knot $Ki$ from the trivial one. Also in \cite{Mih} the
question was placed  whether the group $G_{3}(Ki)$ is able to distinguish the Kishino knot from the trivial one or not. Note that the group $G_{3}(U)$ of the
trivial knot $U$ is isomorphic to $F_2$. In the present section the positive answer on the question will be given. For that purpose the fact that the group
$G_{3}(Ki)$ is not isomorphic to the free group $F_2$ is proved. To do that the Magnus representation of a free metabelian group by matrices $2 \times 2$ is
required. Recall the Magnus representation (see \cite[Chapter 1, Section 4]{LS}).

Let $R$ be a commutative ring over $\mathbb{Z}$ containing independent invertible elements $s_1,$ $\ldots$, $s_q$ and $t_1,$ $\ldots$, $t_q$. Magnus showed that
the map $\varphi$ from the free metabelian group $F_q / F_q''$ into $\mathrm{GL}_2(R)$ defined by $$ x_i \varphi = \left(
  \begin{array}{cc}
    s_i & t_i \\
   0 & 1 \\
  \end{array}
\right), ~~ i = 1, \ldots, q, $$ is a monomorphism.

Using the Magnus representation, the following theorem is stated.

\begin{theorem} The group $G=G_3(Ki)$ having generators $a$, $b$, $c$, $d$ and the system of defining relations $$ d^{-1} l^{-d} c^{-2d^{-1}} b^{-d} c^{-2d^{-1}} a
a^{-2d} d=a^{-1} b^{-d} c^{-2d^{-1}} a, \eqno{(1)} $$ $$ c^{-1} b c =  b^{-d} c^{d^{-1}} b^{d}, \eqno{(2)} $$ $$ c= b^{-d} c^{-2d^{-1}} b^{-d} c^{-2d^{-1}} a
a^{-d} a^{-1}  c^{2d^{-1}}b^{2d}. \eqno{(3)} $$ is not isomorphic to the free group of rank 2. \end{theorem}

\textbf{Proof.} Proof by contradiction. Assume that  $G$ is generated freely by some elements $x$, $y$, i.~e.  $G= \langle x,y \rangle \cong F_2$. Consider the
Magnus representation $\varphi :  G \longrightarrow \mathrm{GL}_2(R)$ defined by the action on generators: $$ x \mapsto \left( \begin{array}{cc}
  s_1 & t_1 \\
  0 & 1 \\
\end{array} \right), \quad y \mapsto \left( \begin{array}{cc}
  s_2 & t_2 \\
  0 & 1 \\
\end{array} \right). $$ Here $s_1$ and $s_2$ are free generators of a free abelian group of rank 2. Consider elements $t_1$, $t_2$ as generators of a free
module $M$ over the Laurent  ring $\mathbb{Z}[s_1^{\pm 1}, s_2^{\pm 1}]$.

The representation maps the generators $a$, $b$, $c$, $d$ of the group $G$ to  matrices $$ a \mapsto \left( \begin{array}{cc}
  A & \gamma \\
  0 & 1 \\
\end{array} \right), \quad b \mapsto \left( \begin{array}{cc}
  B & \lambda \\
  0 & 1 \\
\end{array} \right), \quad c \mapsto \left( \begin{array}{cc}
  C & \mu \\
  0 & 1 \\
\end{array} \right), \quad d \mapsto \left( \begin{array}{cc}
  D & \nu \\
  0 & 1 \\
\end{array} \right), $$ where $A, B, C, D \in \langle s_1, s_2 \rangle$, $\gamma, \lambda, \mu, \nu \in M$.

The projection of $\varphi(G)$ to the entry (1,1) defines the homomorphism of the group $G$ onto the group  $\langle A, B, C, D\rangle \leq \langle s_1, s_2
\rangle$ with the kernel $G'$. Using relations (1)--(3) in the quotient group $G/G'$, the following relations are obtained: $$ ABC^2=1,\quad B=C,\quad AC^3=1.
$$ Thus $$ A=C^{-3},\quad B=C. $$ Therefore, the representation $\varphi$ acts on the generators of the group $G$ in the following manner: $$
 a \mapsto
\left( \begin{array}{cc}
  C^{-3} & \gamma \\
  0 & 1 \\
\end{array} \right), \quad b \mapsto \left( \begin{array}{cc}
  C & \lambda \\
  0 & 1 \\
\end{array} \right), \quad c \mapsto \left( \begin{array}{cc}
  C & \mu \\
  0 & 1 \\
\end{array} \right), \quad d \mapsto \left( \begin{array}{cc}
  D & \nu \\
  0 & 1 \\
\end{array} \right). $$ Is it clear to see that $\langle C,D \rangle = \langle s_1, s_2 \rangle$.

As a further step consider the relation (2). Using the homomorphism $\varphi$, the equality is obtained $$ \left( \begin{array}{cc}
  C & C^{-1}\lambda +(1-C^{-1})\mu \\
  0 & 1 \\
\end{array} \right) = $$ $$ =\left( \begin{array}{cc}
  C & D^{-1}(1-C^{-1})\lambda+ C^{-1}D\mu +(D^{-1}(C-1)+(C^{-1}-1)(D^{-1}+1))\nu \\
  0 & 1 \\
\end{array} \right). $$ Hence, in the module $M$ there is the relation $$ C^{-1}\lambda +(1-C^{-1})\mu= D^{-1}(1-C^{-1})\lambda+ C^{-1}D\mu
+(D^{-1}(C-1)+(C^{-1}-1)(D^{-1}+1))\nu $$ or $$ (C^{-1}+ D^{-1}(C^{-1}-1))\lambda + (1-C^{-1} - C^{-1}D)\mu+ (1-C^{-1})(D^{-1}+1- D^{-1}C)\nu=0. $$ Multiplying
by $CD$, we obtain $$ (D+ 1-C)\lambda + D(C-1 - D)\mu+ (C-1)(1+ D-C)\nu=0, $$ which is equivalent to the relation $$ (1-C+D) \left( \lambda -D\mu+(C-1)\nu
\right)=0. $$ Since the module $M$ is torsion free $$
 \lambda -D\mu+(C-1)\nu =0.
$$ Rewrite this relation in the form $$
 D\mu+\nu= \lambda +C\nu.
 \eqno{(4)}
$$

Acting by the homomorphism $\varphi$ on the products $dc$ and $bd$, we obtain $$ \varphi(dc)= \left( \begin{array}{cc}
  CD & D\mu+\nu \\
  0 & 1 \\
\end{array} \right), \quad \varphi(bd)= \left( \begin{array}{cc}
  CD & \lambda +C\nu \\
  0 & 1 \\
\end{array} \right). $$ Due to (4) we have the equality $\varphi(dc)=\varphi(bd)$, i.~e. $\varphi(c^{-1}b^d)=1$. Since the Magnus representation is a faithful
representation of the free metabelian group, the following relation ought to hold in the quotient group $G/G''$ $$ c =  b^{d}. \eqno{(5)} $$

Let us prove the following lemma.

\begin{lemma} \label{8} The quotient group of the group $G$ by the normal closure of $c^{-1}b^d$ is isomorphic to the group $$ G_1= \left\langle \,a, c, d ~ \|
~ c^{-1} c^{-2d^{-1}} = a a^{d} a^{-1}\, \right\rangle . $$ \end{lemma}

\textbf{Proof.} After the substitution $b^{d}=c$ the relation (2) becomes $$ c^{-1} b c =  c^{-1} c^{d^{-1}} c \Rightarrow
 b  =  c^{d^{-1}} \Leftrightarrow c =  b^{d},
$$ i.~e. it is identity due to (5).

After the substitution $b^{d}=c$ the relation (3) is of the form $$ c= c^{-1} c^{-2d^{-1}} c^{-1} c^{-2d^{-1}} a a^{-d} a^{-1}  c^{2d^{-1}}c^{2} \Rightarrow $$
$$ 1 =  c^{-2d^{-1}} c^{-1} c^{-2d^{-1}} a a^{-d} a^{-1}  c^{2d^{-1}}
 \Leftrightarrow
$$ $$
 c^{-1} c^{-2d^{-1}} = a a^{d} a^{-1}.
 \eqno{(6)}
$$

After the substitution $b^{d}=c$ the relation (1) becomes $$ d^{-1} c^{-1} c^{-2d^{-1}} c^{-1} c^{-2d^{-1}} a a^{-2d} d=a^{-1} c^{-1} c^{-2d^{-1}} a. $$
Substituting the product $c^{-1} c^{-2d^{-1}} = a a^{d} a^{-1}$ due to (6), we obtain $$ d^{-1} a a^{d} a^{-1} a a^{d} a^{-1} a a^{-2d} d=a^{-1} a a^{d} a^{-1}
a \Rightarrow $$ $$ d^{-1} a a^{d} a^{d}  a^{-2d} d= a^{d}  \Rightarrow $$ $$ d^{-1} a  d= a^{d} $$ --- the identity.

Thus the quotient group has only one defining relation (6). \hfill $\square$

\medskip

Denote by $F_3= \langle a, c, d \rangle$ a free group of rank 3 with free generators $a,c,d$. Due to the assumption on the group $G$ and in view of the lemma
\ref{8} the group $$
 \left\langle \,a, c, d ~ \| ~ c c^{2d^{-1}} a a^{d} a^{-1}=1,~ F_3'' \right\rangle
$$ is isomorphic to the free metabelian group $ F_2/F_2''$ of rank 2.

Due to results of the paper  \cite{GGN} it is possible if and only if the element $$ w=c c^{2d^{-1}} a a^{d} a^{-1} \in F_3 $$ is a primitive in the free
metabelian group  $F_3/F_3''$. Due to \cite{Ro} the element $w$ is primitive in $F_3/F_3''$ if and only if the vector $$ (\partial_a w, \partial_c w,\partial_d
w) $$ is unimodular (i.~e. the ideal generated by its components contains 1 in the Laurent  ring $\mathbb{Z}[a^{\pm 1}, c^{\pm 1}, d^{\pm 1}]$). Here
 $\partial_a w, \partial_c w,\partial_d w$ are left Fox derivatives by variables $a, c, d$ respectively (see \cite[Chapter 1, Section 10]{LS}).

Values of Fox derivatives are taken in the Laurent ring $\mathbb{Z}[F_3/F_3']$. Hence, to simplify the notation the assumption $F_3/F_3'= \langle a, c, d \rangle$
is made, i.~e. $a, c, d$ are generators of the free abelian group. We have $$ \partial_a w=c^3(1+d^{-1}a-a), $$ $$ \partial_c w=1+cd+c^2d, $$ $$ \partial_d
w=c-c^3-c^3ad^{-1}+c^3a^2d^{-1}. $$ Pass from the vector $(\partial_a w, \partial_c w,\partial_d w)$ to the vector $$ (c^{-3}d\partial_a w, \partial_c w,
c^{-1}d\partial_d w)= $$ $$ =( a+d-ad, 1+cd+c^2d, d-c^2d-c^2a+c^2a^2 ). $$ Due to assumption on the unimodularity there are  coefficients $$ A, C, D \in
\mathbb{Z}[a^{\pm 1}, c^{\pm 1}, d^{\pm 1}] $$ such that the following equation holds $$ (a+d-ad)A+ (1+cd+c^2d)C+ (d-c^2d-c^2a+c^2a^2)D=1. $$ Multiplying by a
suitable monomial of the form $a^{\gamma}c^{\pi}d^{\nu}$, where $\gamma, \pi, \nu $ are positive integers, we obtain $$ (a+d-ad)A'+ (1+cd+c^2d)C'+
(d-c^2d-c^2a+c^2a^2)D'=a^{\gamma}c^{\pi}d^{\nu}
 \eqno{(7)}
$$ which holds in the polynomial ring, $A',C',D' \in \mathbb{Z}[a, c, d]$.

Consider the system of equations $$ a+d-ad=0,
 \eqno{(8)}
$$ $$ 1+cd+c^2d=0,
 \eqno{(9)}
$$ $$ d-c^2d-c^2a+c^2a^2=0.
 \eqno{(10)}
$$ It follows from (9) that $$ d=- \frac{1}{c^2+c}.
 \eqno{(11)}
$$ From (8) due to (11) we have $$ a=\frac{1}{c^2+c+1}.
 \eqno{(12)}
$$ From (10) due to (11) and (12) the equation on $c$ is obtained $$ c^3-c^2-c-1=0.
 \eqno{(13)}
$$

Let $c_0$ be some solution of (13), generally speaking a complex number. Since system of equations $$ \left\{ \begin{array}{l}
  c^2+c=0, \\
  c^3-c^2-c-1=0 \\
\end{array} \right. \quad \mbox{and} \quad \left\{ \begin{array}{l}
  c^2+c+1=0, \\
  c^3-c^2-c-1=0 \\
\end{array} \right. $$
are not compatible, $$ d_0=- \frac{1}{c_0^2+c_0}, \quad a_0=\frac{1}{c_0^2+c_0+1} $$ are well defined. Note also that $a_0,c_0,d_0 \neq
0$. Substituting $a=a_0,c=c_0,d=d_0$ in (7), we obtain $$ 0=a_0^{\gamma}c_0^{\pi}d_0^{\nu}\neq 0 $$
 which is a contradictory expression. Thus the vector $(\partial_a w, \partial_c w,\partial_d w)$ is not unimodular. Therefore, the group $G$ is not a free
 group. \hfill $\square$

\vspace{0.8cm}

\section{Representations of
groups by formal power series}

\vspace{0.5cm}

In the present section we will show that for the studying of virtual knot groups we can use representations by rings of formal power series with non-commutative variables.

Recall (see, for example \cite[Chapter 1, Proposition 10.1]{LS}) a well known representation of a free group $F_n = \langle x_1, x_2, \ldots, x_n \rangle$ into
the ring of formal power series $\mathbb{Z}[[X_1,$ $X_2,$ $\ldots, X_n]]$ of noncommutative variables $X_1, X_2, \ldots, X_n$ defined by the action on
generators: $$ x_i \mapsto 1 + X_i, ~~i = 1, 2, \ldots, n. $$ In that case inverse elements of generators go to the following formal power series $$ x_i^{-1}
\mapsto 1 - X_i + X_i^2 - X_i^3 + \ldots , ~~i = 1, 2, \ldots, n. $$ The representation defined in that manner is faithful  (i.~e. its kernel is
trivial). Moreover, it remains being faithful if the ring $\mathbb{Z}[[X_1, X_2, \ldots, X_n]]$ is replaced by the quotient ring $\mathbb{Z}[[X_1, X_2,
\ldots, X_n]] / \langle X_1^2, X_2^2, \ldots, X_n^2 \rangle$ by the two-sided ideal generated by elements $X_1^2, X_2^2, \ldots, X_n^2$.

Let $$ \mathcal{P} = \left\langle \, x_1,\ldots,x_n     \,\| \, r_1,\ldots,r_m         \, \right\rangle $$ be some finite presentation of the group $G$ and
$A_n=\mathbb{Q}[[X_1,\ldots,X_n]]$ is an algebra of formal power series of noncommutative variables $X_1,\ldots,X_n$ over the field of rational numbers. Define series $f_j$, $j=1,\ldots,m$, in
the algebra $A_n$ by equalities $$ f_j=r_j(1+X_1,\ldots,1+X_n)-1, \quad j=1,\ldots,m. $$

\smallskip

\begin{prop} The quotient algebra $A_n/ \langle f_1, \ldots, f_m \rangle$ is the invariant of the group $G$, i.~e. it does not depend on the explicit
presentation. \end{prop}

\textbf{Proof.} It is sufficient to show that the quotient algebra  $A_n/ \langle f_1,\ldots,f_m \rangle$ is unaltered in case the presentation $\mathcal{P}$ is
changed by Tietze transformations.

1. In case of bringing in a new generator $z=w(x_1,\ldots,x_n)$
 a new variable $Z$ and the relation of the form $Z=W(X_1,\ldots,X_n)$ are added to the algebra
$$ A_n/ \langle f_1,\ldots,f_m \rangle, $$ where $W(X_1,\ldots,X_n)=w(1+X_1,\ldots,1+X_n)-1$. Hence the quotient algebra $A_n/ \langle f_1,\ldots,f_m \rangle$
is unaltered. Analogously in case of excluding a generator the quotient algebra $A_n/ \langle f_1,\ldots,f_m \rangle$ stays unaltered.

2. In case of bringing in the relation of the form $r_p r_q$ in the group $G$ the relation $f_p +f_q-f_pf_q$ is added to the algebra, which does not change the
ideal $\langle f_1,\ldots,f_m \rangle$. By the analogy in case of excluding the relation the quotient algebra  $A_n/ \langle f_1,\ldots,f_m \rangle$ is
unaltered. \hfill $\square$

\medskip

Considering quotients of $A_n$ by some ideals, we can construct representations of groups in finite dimensional algebras.
For example, consider the quotient   $B_2 = \mathbb{Q}[[X,Y]]/ \langle X^2, Y^2 \rangle$ and prove

\begin{prop} An arbitrary nontrivial quotient algebra of the algebra $B_2$ has a finite dimension.
 \end{prop}

\textbf{Proof.} Let $f=f(X,Y)$ be a nontrivial polynomial in the algebra $B_2$. Then its summand of the highest power has one of the following form $$
a(XY)^k+b(YX)^k, \quad a(XY)^kX+b(YX)^kY, \quad a,\, b \in \mathbb{Q}. $$ Assume at first that $$ f=g+a(XY)^k+b(YX)^k, $$ where the power of element $g \in B_2$
is strictly less than $2k$. Suppose that $a\neq 0$. Multiplying this equation by $X$ on the right and by $Y$ on the left, we obtain $$
fX=gX+a(XY)^kX, \quad Yf=Yg+aY(XY)^k. $$ Therefore, in the quotient algebra $B_2 / \langle f \rangle$ monomials of the power at least $2k+1$ can be
expressed through monomials of the power at most $2k$, i.~e. the dimension of the quotient algebra $B_2/ \langle f \rangle$ is at most $4k+1$.

The second case: 
 $$
  f=g+a(XY)^kX+b(YX)^kY
 $$
 is considered by analogy. \hfill $\square$

\medskip

Let us construct representations of the group $G_{1,r}$ and  $G_{2}$ in $B_2$. Let the generators of $G_{1,r}$ go to the elements:
 $$ x \mapsto 1+X,\quad y \mapsto 1+Y $$  of  $B_2$. Then $$ x^{-1} \mapsto 1-X,\quad y^{-1} \mapsto 1-Y. $$
Find out the additional relation on variables $X$, $Y$ in such way that the map $$ x \mapsto 1+X,\quad y \mapsto 1+Y $$ can be extended to the representation of
the group $G_{1,r}$.

If we go to the algebra than the relation of the group $G_{1,r}$ has the form $$ (1-rX)(1-Y)(1+X)(1+Y)(1+rX)=~~~~~~~~~~~~~~~~~~~~~~~~~~~~~~~~~~~ $$ $$
~~~~~~~~~~~~~~~~~~~~~~~~~~~~~~~~~~~~~~=(1-2Y)(1+rX)(1+Y)(1+X)(1-Y)(1-rX)(1+2Y). $$ 
Opening the brackets and collecting similar terms, the following relation is obtained
$$ (XY)^2-(YX)^2Y+r(XY)^3-r^2(YX)^3-2r^2(YX)^3Y=0. $$ This relation should be added to the algebra $B_2$ to get the representation of the group $G_{1,r}$.
Multiply the last relation by $Y$ on the left $$ (YX)^2Y+r(YX)^3Y=0. $$ Since $1+rXY$ is the invertible element,  $(YX)^2Y=0$.
Multiplying the relation by $X$ on the right, we have $$ (XY)^2X-(YX)^3+r(XY)^3X-2r^2(YX)^4=0. $$ Due to $(YX)^2Y=0$ we obtain $(XY)^2X=0$. Thus the initial
relation becomes $ (XY)^2=0.$

Therefore, the following proposition is proved.

\begin{prop} The group $G_{1,r}$ allows the representation $$ x \mapsto 1+X,\quad y \mapsto 1+Y $$ into the group of invertible elements of the
algebra $B_2 / \langle (XY)^2 \rangle$. \end{prop}

\begin{remark}
Note that the algebra $B_2 / \langle (XY)^2 \rangle$ is a vector space over the field $\mathbb{Q}$ of dimension 8 with the basis $1$, $X$, $Y$, $XY$,
$YX$, $XYX$, $YXY$, $(YX)^2$. Hence, it is possible to construct a linear representation of $G_{1,r}$.
\end{remark}

Consider now the group $G_2$ and construct a map into the algebra  $B_2$. Find out the additional relation on variables $X$, $Y$ in such a
manner that the map $$ x \mapsto 1+X,\quad y \mapsto 1+Y $$ can be extended to the representation of the group $G_2$.

The relation of the group $G_2$ passing into the algebra becomes $$ (1+X)(1-Y)(1-X)(1+Y)(1+X)=~~~~~~~~~~~~~~~~~~~~~~~~~~~~~~~~~~~~~~~~~~~~ $$ $$
~~~~~~~~~~~~~~~~~~~~~~~~~~~~~~~~~~~~~~~~=(1-2Y)(1+X)(1+Y)(1-X)(1-Y)(1+X)(1+2Y). $$ Opening brackets and collecting similar  terms, the following relation is obtained
$$ 2XY-2YX-4YXY-(YX)^3+(XY)^3-2(YX)^3Y=0. $$ This relation should be added to the algebra $B_2$ to get the representation of the
group $G_2$.

Multiply the relation by $Y$ on the left $$ 2YXY+(YX)^3Y=0. $$ Since $2+(XY)^2$ is the invertible element, $YXY=0$.
Multiplying the relation by $X$ on the right, we have $$ 2XYX-4(YX)^2+(XY)^3X-2(YX)^4=0. $$ Due to the equality $YXY=0$ we obtain $XYX=0$. Hence the initial
relation becomes of the form $ XY=YX. $

Therefore, it has been  proven

\begin{prop} The group $G_2$ allows the representation $$ x \mapsto 1+X,\quad y \mapsto 1+Y $$ into the algebra  $B_2 / \langle XY-YX \rangle$.
\end{prop}

Note that the algebra $\mathcal{A}_2 / \langle XY-YX \rangle$ is a vector space over the field $\mathbb{Q}$ of dimension 4 with the basis $1$, $X$, $Y$, $XY$.

\vspace{0.8cm}
\section{Open questions }

\vspace{0.5cm}

In conclusion we formulate series of open questions.

If $K$ is a classical knot than $G_{1,1}(K)$ is a classical knot group, i.~e. a fundamental group of the compliment of $K$ in 3-dimensional sphere. Nelson and
Neumann proved \cite{NN} that the group $G_{1,1}(K)$ is a subgroup of the group $G_{1,r}(K)$ for arbitrary $r \geq 1$ and in contrast to $G_{1,1}(K)$ it
determines the knot $K$ up to reflection.


\begin{problem} 1) Is it true for links?

2) Is it possible to prove the similar statement for groups $G_{3}(K)$, which constructed by using the extension of the Wada representation $W_{3}$?

\end{problem}

Let $G_2$ be the virtual trefoil group defined above.

\begin{problem}

1) Does the equation $\gamma_{\omega} G_2=\gamma_{\omega+1} G_2$ hold?

2) Is there the isomorphism $$ \gamma_k G_2/ \gamma_{k+1} G_2 \cong (\gamma_k F_2/ \gamma_{k+1} F_2 ) \otimes \mathbb{Z}_4 $$ for arbitrary integer $k \geq 2$?
\end{problem}

\begin{problem}
Is there a virtual knot $K$ for which some group from the set $\{ G_{1,1}(K), G_{2}(K), G_{3}(K) \}$ is residually nilpotent?
\end{problem}

\end{document}